\documentclass{elsart}
\usepackage{amsfonts}
\newtheorem{theorem}{Theorem}

\newtheorem{proposition}{Proposition}
\newtheorem{remark}{Remark}
\newcommand{\pp}{\noindent {\em Proof. }}
\newcommand{\bee}[1]{\begin{equation}\label{#1}}
\newcommand{\beq}[1]{\begin{eqnarray}\label{#1}}
\newcommand{\ene}{\end{equation}}
\newcommand{\eqe}{\end{eqnarray}}
\newcommand{\ben}{\begin{eqnarray*}}
\newcommand{\eqn}{\end{eqnarray*}}
\newcommand{\mt}[4]{\left[\begin{array}{cc}#1&#2\\#3&#4\end{array}\right]}

\newcommand{\tr}[1]{ ^t\!#1}

\newcommand{\hg}{\widehat{\Gamma}}
\newcommand{\G}{\Gamma}
\newcommand{\g}{\gamma}

\newcommand{\vp}{\varphi}

\newcommand{\iv}[1]{#1^{-1}\!}
\newcommand{\ve}{\varepsilon}

\newcommand{\bc}{\bar{\chi}}

\newcommand{\su}[1]{\mathrm{Supp}\,#1}
\newcommand{\Aut}[1]{\mathrm{Aut}\,#1}
\newcommand{\Sp}[1]{\mathrm{Span}\,\{#1\}}

\begin{document}

\begin{frontmatter}



\title{Group Gradings on $G_2$.}


\author[yb]{Y. A. Bahturin}
\address{Department of Mathematics and Statistics\\Memorial
University of Newfoundland\\ St. John's, NL, A1C5S7,
Canada}
\ead{yuri@math.mun.ca}\thanks[yb]{Work is partially supported by NSERC grant 227060-04}

\author[mt]{M. V. Tvalavadze}
\address{Department of Mathematics and Statistics\\Memorial
University of Newfoundland\\ St. John's, NL, A1C5S7, Canada}\ead{marina@math.mun.ca}
\thanks[mt]{Work is partially supported by Atlantic Algebra Centre and NSERC grants}

\begin{abstract}
In this paper we describe all group gradings by a finite abelian group
$\Gamma$ of a simple Lie algebra of type $G_2$ over an
algebraically closed field $F$ of characteristic 0.
\end{abstract}

\end{frontmatter}

\section{Introduction}\label{s1}

In this paper we describe, up to an isomorphism, all group gradings on a simple Lie
algebra $L$ of the type $G_2$ over an algebraically closed field
$F$ of characteristic zero. It is well known
(see, for example, \cite{PZ}) that the elements of the support
of the grading in the case of a simple Lie algebra must commute, and so one may always assume that the grading group
$\Gamma$ is finitely generated abelian. If $\Gamma$ is torsion free then it is well-known that a Cartan subalgebra $H$ of $L$ is contained in the identity component of the grading. In this case any root
space is $\Gamma$-graded so that any graded component is the sum
of the root subspaces. On the other hand, gradings by finite cyclic groups have been completely described by Victor Kac \cite{VK}(see also a very useful book
\cite{VO}). As mentioned in \cite{VO}, the so called \textit{Jordan gradings} of $G_2$ have been determined in 1974 by A. Alekseevsky \cite{A}.

In the case of arbitrary gradings by abelian groups, an important reference  is the paper \cite{AE} where the author determines the equivalence classes of abelian group gradings of the octonions. In this paper we use the realization of $L$ as the Lie algebra of derivations of the split octonions \cite{JLA} to obtain all possible gradings of $G_2$ by $\Gamma$ when this group is \emph{finite abelian}. An independent paper is \cite{DM} where the authors also produce a list of arbitrary gradings on $G_2$. Both papers use \cite{AE}, but in a rather different way, maybe thanks to a different approach to the definition of a grading in general.

\section{Basic Facts and Notation}\label{s2}

Given an algebra $A$ over a field $F$, and a group $\G$, we say
that $A$ is $\G$-graded if $A=\bigoplus_{\gamma\in \G}A_{\gamma}$
where each $A_{\gamma}$ is a vector subspace of $A$ and $A_{\gamma}A_{\delta}\subset A_{\gamma\delta}$,
for any $\gamma,\delta\in\G$. A subspace (subalgebra, ideal) $B$ of $A$ is called \textit{graded}
if $B=\bigoplus_{\gamma\in \G}(B\cap A_{\gamma})$. The set $\su{A}=\{\gamma\in\G\,|\, A_{\g}\neq \{ 0\}\}$
is called the \textit{support} of the above grading. For any abelian group $\G$ we denote by $\hg$
the group of multiplicative characters of $\G$, that is, the group homomorphisms $\chi:\G\rightarrow F^{\ast}$.
If $F$ is algebraically closed of characteristic zero and $\G$ is finite abelian, then $\G$-gradings are completely determined by $\hg$-actions in the sense that a subspace $B$ of $A$ is $\G$-graded
if and only if $B$ is invariant under the action of $\hg$ by the automorphisms of $A$. The action in the case
 of a $\G$-grading is defined by setting
\begin{equation}\label{eaction}
\chi\circ\sum_{\g\in\G}a_{\g}=\sum_{\g\in\G}\chi(\g)a_{\g}\mbox{ where }a_{\g}\in A_{\g}\mbox{ for each } \g\in\G.
\end{equation}
Conversely, if there is an action of $\hg$ by automorphisms, we set
$$
A_\g=\{ a\in A\,|\, \chi\circ a=\chi(\g)a\mbox{ for all }\g\in\hg\}.
$$
Note that the natural embedding of $\G$ in $\Aut A$ defined by this action is injective if and only if $\Gamma$ is generated by $\su{A}$. We will always assume that this condition holds for our gradings.

Two $\G$-gradings $A=\oplus_{g\in
\G} A_g$ and $A=\oplus_{\g\in \G} \widetilde{A}_{\g}$ are called {\it
isomorphic} if there exists an automorphism $\sigma: \G\to \G$ and
an automorphism $f: A\to A$ such that $f(A_\g)=\widetilde{A}_{\sigma(\g)}$ for all
$\g\in \G$.

We recall that the split octonion algebra over a field $F$ can be interpreted as the set of matrices
\begin{equation}\label{eoct}
C=\left\{\left.\left[\begin{array}{cc}\alpha&u\\v&\beta\end{array}\right]\right|\alpha,\beta\in F, u,v\in F^3\right\},
\end{equation}
with the product \bee{eocpr}
\mt{\alpha}{u}{v}{\beta}\mt{\alpha'}{u'}{v'}{\beta'}=\mt{\alpha\alpha'-(u,v')}{\alpha
u'+\beta'u+v\times v'}{\alpha'v+\beta v'+u\times
u'}{\beta\beta'-(v,u')} \ene where $u\times v$ is the standard
cross product in $F^3$, and $(u,v)$ is the standard inner product
(see \cite{JLA}). The derivation algebra $L=\mathrm{Der}\, C$ is
the simple Lie algebra of type $G_2$. Let us set
$e_1=\mt{1\;}{0\;}{0\;}{0\;}$, $e_2=\mt{0\;}{0\;}{0\;}{1\;}$,
$u_i=\mt{0\;}{\ve_i}{0\;}{0\;}$, and
$v_j=\mt{0\;}{0}{\ve_j}{0\;}$, where $1\le i,j\le 3$. Here
$\{\ve_1,\ve_2,\ve_3\}$ is the standard basis in $F^3$ with the
cross product given by $\ve_1\times\ve_2=-\ve_2\times\ve_1=\ve_3$,
$\ve_2\times\ve_3=-\ve_3\times\ve_2=\ve_1$ and
$\ve_3\times\ve_1=-\ve_1\times\ve_3=\ve_2$ and the inner product
$(\ve_i,\ve_j)=\delta_{ij}$, the Kronecker delta. The basis
$\{e_1,e_2,u_1,u_2,u_3,v_1,v_2,v_3\}$ is called {\it standard}
with the multiplication table given by:
$$
\begin{array}{c|cccccccc} \hline
     &  e_1  & e_2   & u_1  & u_2  & u_3 & v_1 & v_2 &  v_3\\
\hline
 e_1 &  e_1  & 0     & u_1  & u_2  & u_3 &   0 &  0  & 0\\
 e_2 &  0    & e_2   & 0    &  0   & 0   & v_1 & v_2 & v_3\\
u_1  &  0    & u_1   & 0    & v_3  & -v_2& e_1 &  0   &0\\
u_2  &  0    & u_2  & -v_3 &  0   &  v_1&  0  &  e_1 & 0\\
u_3  &  0    &  u_3  & v_2  & -v_1 &   0  & 0  &   0 & e_1\\
v_1  &  v_1  &  0    &  e_2 & 0    &   0  & 0  & -u_3 & u_2\\
v_2  &  v_2  & 0     & 0    & e_2  &   0  & u_3&  0   & -u_1\\
v_3  & v_3   & 0     & 0    & 0    &   e_2& -u_2& u_1 &0
\end{array}
$$
Note that $1=e_1+e_2$ is the identity element in $C$.

A canonical basis of $L=\mathrm{Der}\, C$ can then be given in the
following way. We recall that for any $x,y\in C$ the mapping
$D_{x,y}=[x_L,y_L]+[x_L,y_R]+[x_R,y_R]$, where $x_L$, $x_R$ are
the left and right multiplications by $x$ in $C$, is a derivation
of $C$ called the \textit{inner} derivation. It is well-known that
all derivations of $C$   are inner. If $T$ is a traceless $3\times
3$-matrix then the mapping \bee{edt}
d_T:\mt{\alpha}{u}{v}{\beta}\rightarrow \mt{0}{uT}{-vT^t}{0}
\ene is also a derivation of $C$ and the set of all such
derivations forms a subalgebra $S$ of $L$ isomorphic to
$\mathrm{sl}\,(3)$. It is claimed in \cite[p.143]{JLA} that the
derivations $D_{e_1,u}$, $D_{e_2,v}$, $d_T$, where $u\in
U=\Sp{u_1,u_2,u_3}$, $v\in V=\Sp{v_1,v_2,v_3}$ and $T\in
\mathrm{sl}\,(3)$ span the whole of $L$. So we can graphically view
$L$ as the set of $4\times 4$-matrices of the form
$\mt{T}{\tr{u}}{v}{0}$.

Now it is known from \cite[p. 285]{JLA} that any automorphism of
$L$ can be written in the form $D\rightarrow ADA^{-1}$ where $A$
is an automorphism of $C$. Thus the following is true.

\begin{proposition}\label{p1}
For any grading $L=\bigoplus_{\g\in\G}L_\g$ there is exactly one grading
$C=\bigoplus_{\g\in\G}C_\g$ such that $L_{\delta}=L\cap (\mathrm{End}\,C)_{\delta}$
for any $\delta\in\G$. Here
$$(\mathrm{End}\,C)_{\delta}=\{\vp:C\rightarrow C\,| \mbox{ such that } \vp(C_\g)\subset C_{\delta\g}\mbox{ for any }\g\in\G\}.
$$
Two $\Gamma$-gradings of $L$ are isomorphic if and only if there corresponding gradings of $C$ are isomorphic.
\end{proposition}

\pp It is obvious that if $L_{\delta}$ is defined in the way as claimed,
then we obtain a $\G$-grading on $L$.

Now let us assume that we have a $\G$-grading on $L$. Then $\hg$ acts by
the automorphisms $\bc$ for each $\chi\in\hg$, defined by $\bar{\chi}(x)=\chi\circ x$ (see (\ref{eaction})).
To each $\bc$ there is an automorphism $A_{\chi}:C\rightarrow C$ such that $\bc(D)=A_{\chi}D\iv{A_{\chi}}$.
Such $A_{\chi}$ is defined uniquely because according to \cite[p.43]{JLA} there is a subspace
 $C_0$ in $C$ of codimension 1, which is invariant under $L$. Now if $A_{\chi}$ induces a
 trivial mapping on $L$ then by Schur's Lemma $A_{\chi}$ is scalar on $C_0$
 with coefficient $\lambda$. Recall that $C_0$ is the subspace of the octonions $\mt{\alpha}{u}{v}{\beta}$ with $\alpha+\beta=0$.
 Now $A_{\chi}(u_1u_2)=A_{\chi}(v_3)$, $A_{\chi}(u_1)A_{\chi}(u_2)=A_{\chi}(v_3)$ implies $\lambda^2=\lambda$.
 So we must have $\lambda=1$. In this case the mapping  $\chi\mapsto A_{\chi}\in\mathrm{Aut}\, C$ is a group homomorphism.
It follows that there is a grading by $\G$ on $C$. But $D\in L_{\gamma}$ if and only if $\bc(D)=\chi(\g)D$
 for any $\chi\in\G$. Now $\bc(D)=A_{\chi}D\iv{A_{\chi}}$. If $z\in C_{\g}$ then $A_\chi(z)=\chi(\gamma)z$ for any $\chi\in\hg$.
 Thus
$$
A_{\chi}(D(z))=(A_{\chi}D\iv{A_{\chi}})(A_{\chi}(z))=\bc(D)(A_{\chi}(z))=\chi(\delta)\chi(\g)D(z)=\chi(\delta\g)D(z)
$$
and so $D(z)$ has degree $\delta\gamma$.

Now if two $\G$-gradings of $C$ are isomorphic by means of an automorphism $A\in \Aut C$ and $\sigma\in \Aut \G$ the obviously the respective gradings of $L$ are isomorphic by $f: D\mapsto ADA^{-1}$, $D\in L$, and the same $\sigma$. Actually, as shown just above, since an isomorphism $f$ of $L$ uniquely defines $A_f\in\Aut C$, the converse is also true. \hfill$\Box$

As a result, if we know all group gradings on the octonions, up to isomorphism, then using Proposition \ref{p1},
we can completely describe all group gradings on $L$, up to isomorphism.

\section{Gradings on octonions}\label{s3}

In \cite{AE} the author obtains a full description of
equivalence classes of gradings on $C$. In that paper two gradings $A=\oplus_{g\in \G} A_g=\oplus_{h\in \Lambda} \widetilde{A}_h$ by finite
groups $\G$ and $\Lambda$ of an algebra $A$ are said to be {\it
equivalent} if there is an automorphism $\varphi$ of $A$ such that
for any $g\in \G$ with $A_g\ne 0$ there is an $h\in \Lambda$ with
$\varphi(A_g)=\widetilde{A}_h$. It may happen that two equivalent
gradings correspond to non-isomorphic groups $\Lambda$ and $\G$.
In other words, the same decomposition as a direct sum of
subspaces may give gradings by non-isomorphic groups.

In this paper we are mainly interested in the description of group
gradings on $G_2$ up to an isomorphism. In order to apply Proposition 1, we need to know the
description of all gradings on the octonions up to an isomorphism of this kind. The
following restatement of Theorem 8 from \cite{AE} provides us with
a full description of group gradings on $C$ up to
isomorphism.

\begin{theorem}
Let $\G$ be a finite Abelian group, and $C=\oplus_{g\in \G} C_g$
be a $\G$-grading of the split octonion algebra $C$. Then, as a $\G$-graded
algebra, $C$ is isomorphic to one of the following:
\par\medskip
\emph{Type 1.} $C=C_e\oplus C_g\oplus C_{g^{-1}}\oplus C_h\oplus
C_{gh}\oplus C_{g^{-1}h}$ where $h\in \G$ is an arbitrary element
of order 2, $g\in \G$ is an arbitrary element of order $>2$, and
all elements $\{e,g,g^{-1},h,gh,g^{-1}h\}$ are different.
Moreover, $C_e=Fe_1+Fe_2$, $C_g=Fu_1$, $C_{g^{-1}}=Fv_1$,
$C_h=Fu_3+Fv_3$, $C_{gh}=Fv_2$, $C_{g^{-1}h}=Fu_2.$
\par\medskip
\emph{Type 2.} $C=C_e\oplus C_g\oplus C_h\oplus C_{gh}\oplus
C_{g^{-1}}\oplus C_{h^{-1}}\oplus C_{g^{-1}h^{-1}}$ where $g,h\in
\G$ are of order $>2$, and all elements
$\{e,g,h,gh,g^{-1},h^{-1},g^{-1}h^{-1}\}$ are different. Moreover,
$C_e=Fe_1+Fe_2$, $C_g=Fu_1$, $C_h=Fu_2$, $C_{gh}=Fv_3$,
$C_{g^{-1}}=Fv_1$, $C_{h^{-1}}=Fv_2$, $C_{g^{-1}h^{-1}}=Fu_3$.
\par\medskip
\emph{Type 3.} $C=C_e\oplus C_h\oplus C_{h^{-1}}\oplus C_{h^2}\oplus
C_{h^{-2}}$ where $h\in \G$ is an arbitrary element of order $>4$.
Moreover, $C_e=Fe_1+Fe_2$, $C_h=Fu_2+Fu_3$,
$C_{h^{-1}}=Fv_2+Fv_3$, $C_{h^{-2}}=Fu_1$, $C_{h^2}=Fv_1$.
\par\medskip
\emph{Type 4.} $C=C_e\oplus C_g\oplus C_{g^{-1}}$ where $g\in \G$
is an arbitrary element of order $>2$. Moreover,
$C_e=Fe_1+Fe_2+Fu_1+Fv_1$, $C_g=Fu_2+Fv_3$,
$C_{g^{-1}}=Fu_3+Fv_2$.
\par\medskip
\emph{Type 5.} $C=C_e\oplus C_g\oplus C_{g^{-1}}$ where $g\in \G$
is an arbitrary element of order 3, $\G=\mathbb{Z}_3$. Moreover, $C_e=Fe_1+Fe_2$,
$C_g=Fu_1+Fu_2+Fu_3$, $C_{g^{-1}}=Fv_1+Fv_2+Fv_3$.
\par\medskip
\emph{Type 6.} $C=C_e\oplus C_g\oplus C_{g^{-1}}\oplus C_{g^2}$ where
$g\in \G$ is an arbitrary element of order 4, $\G=\mathbb{Z}_4$. Moreover,
$C_e=Fe_1+Fe_2$, $C_g=Fu_1+Fu_2$, $C_{g^{-1}}=Fv_1+Fv_2$,
$C_{g^2}=Fu_3+Fv_3$.
\par\medskip
\emph{Type 7. }$C=C_e\oplus C_g$ where $g\in \G$ is an arbitrary element
of order 2, $\G=\mathbb{Z}_2$. Moreover, $C_e=Fe_1+Fe_2+Fu_1+Fv_1$,
$C_g=Fu_2+Fv_2+Fu_3+Fv_3$.
\par\medskip
\emph{Type 8.} $C=C_e\oplus C_g\oplus C_h\oplus C_{gh}$ where $g, h\in \G$, $g\neq h$
are arbitrary elements of order 2, $\G=\mathbb{Z}_2\times \mathbb{Z}_2$. Moreover, $C_e=Fe_1+Fe_2$,
$C_g=Fu_1+Fv_1$, $C_h=Fu_2+Fv_2$, $C_{gh}=Fu_3+Fv_3$.
\par\medskip
\emph{Type 9.} $C=C_e\oplus C_h\oplus C_g\oplus C_k\oplus C_{gh}\oplus
C_{hk}\oplus C_{gk}\oplus C_{hgk}$ where $g,h,k\in \G$ are
arbitrary elements of order 2, $\G=\mathbb{Z}_2\times \mathbb{Z}_2\times \mathbb{Z}_2$.  Moreover, $C_e=F1$,
$C_h=F(e_2-e_1)$, $C_g=F(v_1-u_1)$, $C_k=F(v_2-u_2)$,
$C_{gh}=F(u_1+v_1)$, $C_{hk}=F(u_2+v_2)$, $C_{gk}=F(u_3+v_3)$,
$C_{hgk}=F(v_3-u_3)$.
\end{theorem}
\pp In the case when $\G$ contains at least one element of order
strictly greater than 2, the proof is completely the same as the
proof of Theorem 8 from \cite{AE}. Therefore, we can assume that
there are no elements in $\G$ of order greater than 2. According
to \cite{AE}, the following cases may occur.

{\it Case 1.} Let $\su{\G}$ be generated by one element of
order 2. We denote this element by $g$.  It follows from \cite{AE}
that the grading of $C$ by $\G$ takes the following form: $C_e=H$,
a quaternion subalgebra, and $C_g=Hx$ where $x\in H^{\bot}$ such
$n(x)\ne 0$ where $n$ is the quadratic form on $C$ (see
\cite{NJ}). Next we want to show that all such gradings are indeed
isomorphic. Since $n(x)\ne 0$, normalizing $x$ (if necessary), we
can actually assume that $n(x)=-1$. Further, we choose a basis of
$H$ of the form $\{1,v,w,vw\}$ such that $v^2=w^2=1$. Since $x\in
H^{\bot}$, we have that $x\in (F1)^{\bot}$, and, as a consequence,
$\bar x=-x$ and $x^2=1$. Moreover, for any $h\in H$, $n(x,h)=0$
where $n(x,h)=n(x+h)-n(x)-n(h)$. Due to the relation $n(x)=x\bar
x$, the latter is equivalent to $h\bar x=-x\bar h$. If $t$ denote
any element from the set $\{v,w,vw\}$, then $\bar t=-t$.
Therefore, $tx=-xt$ whenever $t$ is in $\{v,w,vw\}$. Since
$C=H\oplus Hx$, we have that $\{1,v,w,vw,x,vx,wx,(vw)x\}$ is a
basis of $C$. Next we look at the multiplication table of this
basis. Namely,
$$ (vx)(wx)=-(xv)(wx)=-x(vw)x=(vw)x^2=vw,$$
$$ (vx)(vx)=-(xv)(vx)=-x(vv)x=-1,$$
$$ (wx)(wx)=-(xw)(wx)=-x(ww)x=-1,$$
$$((vw)x)(vx)=-(x(vw))(vx)=-x((vw)v)x=(vw)v=-w.$$
In the same way one  can easily find all the remaining products. As
a consequence, the multiplication table for this basis remains the
same no matter what $x\in H^{\bot}$ with $n(x)\ne 0$ is initially
fixed. Finally, if $C=H\oplus Hx=H'\oplus H'x'$ where $H$, $H'$
are two quaternion subalgebras and $x\in H^{\bot}$, $x'\in H'^{\bot}$
with $n(x)\ne 0$, $n(x')\ne 0$, then we can fix bases defined
above for both $H$ and $H'$. These bases can be extended to the bases
of $C$ with the same multiplication tables. Clearly, the mapping
that sends each element of the first basis to the corresponding
element of the second basis is in fact a graded automorphism of
$C$, therefore, the first grading is isomorphic to the second. In
particular, any such grading is isomorphic to the grading of type
7.

{\it Case 2.} Let $\su{\G}$ be generated by two elements of
order 2. We denote these elements by $g$ and $h$, respectively. It
follows from \cite{AE} that the grading of $C$ by $\G$ takes the
following form: $C_e=K$, a two-dimensional composition algebra,
$C_g=Kx$, $C_h=Ky$, $C_{gh}=K(xy)$ where $x\in K^{\bot}$, $y\in
(K+Kx)^{\bot}$ with $n(x)\ne 0$, $n(y)\ne 0$.
 Using the same strategy as in case 1, we can fix a basis of $K$
 of the form: $\{1, v\}$ where $v^2=1$, and, then, extend this basis
 to a basis of $C$ as follows: $\{1, v, x, vx, y, vy, xy, (vx)y\}$. For the same reason
 as above, the multiplication table of this basis does not depend
 on the choice of $x$ and $y$. Any two gradings in this case are
 isomorphic, and, in particular, can be reduced to the grading of
 type 8.

 {\it Case 3.} Let $\su{\G}$ be generated by three elements of
order 2. We denote these elements by $g$, $h$ and $k$,
respectively. It follows from \cite{AE} that the grading of $C$ by
$\G$ takes the following form: $C_e=F1$, $C_h=Fx$, $C_g=Fy$,
$C_k=Fz$, $C_{gh}=F(xy)$, $C_{hk}=F(xz)$, $C_{gk}=F(yz)$,
$C_{hgk}=F(xy)z$ where $x\in (F1)^{\bot}$, $y\in (F1+Fx)^{\bot}$
and $z\in (F1+Fx+Fy+Fxy)^{\bot}$ with $n(x)\ne 0$, $n(y)\ne 0$ and
$n(z)\ne 0$. In the same way as above, we can show that the
multiplication table that corresponds to the basis of $C$ of the
form: $\{1,x,y,z,xy,yz,xz,(xy)z\}$ does not depend on the choice
of $x$, $y$ and $z$. Any two gradings in this case are also
isomorphic, and, in particular, can be reduced to the grading of
type 9. The proof is complete.
 \hfill$\Box$

\section{Elementary Gradings of $G_2$}\label{s4}

Almost all gradings of $C$, as described in Theorem 1, are such
that the subspace $\Sp{e_1,e_2}$ is in the identity component. In
all these cases the subspaces $U$ and $V$ of $C$  are also graded. Actually, each $u_i,\, v_j$ is homogeneous and $\deg u_i=(\deg v_i)^{-1}$, for any $i=1,2,3$.
Then, obviously, the subspaces $D_{e_1,U}$ and $D_{e_2,V}$ are
graded in a compatible way, that is, say, $\deg D_{e_1,u}=\deg u$ for $u$ homogeneous.
Let $(\g_1,\g_2,\g_3)$ be the tuple of degrees of elements $u_1,u_2,u_3$. Then one can easily check that $d_{E_{ij}}\in L_{\g_i^{-1}\g_j}$, for all $i\neq j$ and $d_T\in L_e$ if $T$ is diagonal. So we have a grading on $\mathrm{sl}(3)$ isomorphic to a restriction of an elementary grading of $M_3\cong U\otimes U^{\ast}$ defined by the $3$-tuple $(\g_1,\g_2,\g_3)$ (see \cite{BSZ}). So the tuple $(\g_1,\g_2,\g_3)$
completely defines the $\G$-grading of $L$ in the case of
elementary gradings. Based on Theorem 1, we are able to list all
gradings on $L$ induced by $\G$-gradings on the octonions of Types 1-8.

\subsection{Gradings of {\it Type 1}}\label{ssa}
In this case $U_g=Fu_1$, $U_{g^{-1}h}=Fu_2$, $U_h=Fu_3$. Thus we
have an elementary $\G$-grading on $L$  given by the tuple
$(g,g^{-1}h,h)$, where  $h\in \G$ is an arbitrary element of order
2, $g\in \G$ is an arbitrary element of order $>2$, and all
elements $\{e,g,g^{-1},h,gh,g^{-1}h\}$ are different.

\subsection{Gradings of {\it Type 2} }\label{ssb}
In this case, $U_g=Fu_1$, $U_h=Fu_2$, $U_{g^{-1}h^{-1}}=Fu_3$.
Then the corresponding elementary grading on $L$ is given by the
tuple $(g,h,g^{-1}h^{-1})$, where $g,h\in \G$ are of order $>2$,
and all elements $\{e,g,h,gh,g^{-1},h^{-1},g^{-1}h^{-1}\}$ are
different.

\subsection{Gradings of {\it Type 3}}\label{ssc1}
We have that  $U_h=Fu_2+Fu_3$, $U_{h^{-2}}=Fu_1$. Then an
elementary $\G$-grading on $L$ is defined by $(h^{-2},h,h)$, where $h\in \G$ is an arbitrary element of order $>4$.

\subsection{Gradings of {\it Type 4}}\label{ss3}

In this case, $U_e=Fu_1$, $U_g=Fu_2$, $U_{g^{-1}}=Fu_3$. This
induces an elementary $\G$-grading defined by the tuple
$(e,g,g^{-1})$, where $g\in \G$ is an arbitrary element of order
$>2$.

\subsection{Gradings of {\it Type 5}}\label{ssc2}

In this case, $U_g=Fu_1+Fu_2+Fu_3$. This induces an elementary
grading defined by the tuple $(g,g,g)$, where $g\in \G$ is an
arbitrary element of order 3.

\subsection{Gradings of {\it Type 6}}\label{ssc4}
We have that $U_g=Fu_1+Fu_2$, $U_{g^2}=Fu_3$. Then an  elementary
$\G$-grading on $L$ is given by $(g,g,g^2)$, where
$g\in G$ is an arbitrary element of order 4.

\subsection{Gradings of {\it Type 7}}\label{ssc5}
Here,  $U_e=Fu_1$, $U_g=Fu_2+Fu_3$. Then, an elementary
$\G$-grading on $L$ is given by $(e,g,g)$, where
$g\in G$ is an arbitrary element of order 2.

\subsection{Gradings of {\it Type 8}}\label{ssc6}
In this case, $U_g=Fu_1$, $U_h=Fu_2$, $U_{gh}=Fu_3$. An elementary
$\G$-grading on $L$ is given by $(g,h,gh)$, where $g, h\in \G$, $g\neq h$
are arbitrary elements of order 2.

\section{Non-elementary gradings of $G_2$}\label{s5}

In the case of a $\G$-grading of the Type 9, the grading on $C$ is
obtained by applying the ``Cayley - Dickson'' process to the
composition algebra $F$.
 As a result, up to isomorphism,
 any grading is isomorphic to the following:

 \subsection{Gradings of {\it Type 9}}\label{ssc9}

For this grading all components of $C$ are one-dimensional. None
of $U$, $V$, or $S\cong \mathrm{sl}\,(3)$ are $\G$-graded. An
immediate verification says that $L_e=\{0\}$ and the other
components of the $\G$-grading of  $L$ are given by
\begin{eqnarray*}
&&L_h=\Sp{d_{E_{11}-E_{22}},d_{E_{11}-E_{33}}},\\
&&L_g=\Sp{-D_{e_{1},u_{1}}+D_{e_{2},v_{1}},d_{E_{23}-E_{32}}},\\
&&L_k=\Sp{-D_{e_{1},u_{2}}+D_{e_{2},v_{2}},d_{E_{13}-E_{31}}},\\
&&L_{hg}=\Sp{D_{e_{1},u_{1}}+D_{e_{2},v_{1}},d_{E_{23}+E_{32}}},\\
&&L_{hk}=\Sp{D_{e_{1},u_{2}}+D_{e_{2},v_{2}},d_{E_{13}+E_{31}}},\\
&&L_{gk}=\Sp{D_{e_{1},u_{3}}+D_{e_{2},v_{3}},d_{E_{12}+E_{21}}},\\
&&L_{hgk}=\Sp{-D_{e_{1},u_{3}}+D_{e_{2},v_{3}},d_{E_{12}-E_{21}}}.
\end{eqnarray*}
\begin{theorem}\label{TM} Let $\G$ be a finite abelian group. Any grading by $\G$ on a simple Lie algebra $L$ of type $G_2$ is isomorphic to one of the gradings of the types 1 to 9. Any two gradings belonging to different types 1 to 9 are not isomorphic.
\end{theorem}

We conclude with a couple of remarks.
\begin{remark}\label{r1}
Depending on particular values of the elements of $\G$ satisfying the conditions in the definition of gradings of types 1 to 8 the number and the dimensions of homogeneous components can be different. This explains why the list in \emph{\cite{DM}} is longer than ours. If an elementary grading is given by a tuple $(\g_1,\g_2,\g_3)$ then  $\g^{-1}_i\g_j\in\su{L}$ for any $i,j=1,2,3$. Some of these elements can be equal, which may reduce the number of homogeneous components while increasing their dimensions. But all these gradings are obtained by a clear uniform procedure, so we decided to list them as a single item.
\end{remark}

\begin{remark}\label{r2} If $\G$ is a finite abelian group then a $\G$-symmetric space is a homogeneous space of the form $G/H$ where $G$ is a simply connected Lie group with $\G$ as a subgroup of $\Aut G$ such that $(G^\G)_1\subset H\subset G^\G$. Here $G^\G$ is the set of fixed points of $G$ under the action of $\G$ and $(G^\G)_1$ the connected component of the identity element of $G^\G$. Any such structure leads to the grading of the Lie algebra $\mathfrak{g}$ of $G$ by (a group dual to) $\G$. The converse is also true. Thus, for the classification of such spaces we need to know the structure of gradings on $\mathfrak{g}$ by a particular group $\G$, the point we adopt in this paper.
\end{remark}

\begin{center}
\textbf{Acknowledgment}
\end{center}
The authors thank Professor Efim Zelmanov for his interest in this research and valuable advice. The first author thanks Department of Mathematics of University of California - San Diego for hospitality during his visits as a guest of Professor Efim Zelmanov in 2004 and 2006. The second authors thanks Atlantic Algebra Centre for supporting her research in 2006/07 academic year.

\end{document}